\documentclass{amsart}
\usepackage{amsmath, amsaddr}
\usepackage{amssymb}
\usepackage{amsthm}
\usepackage{graphicx} 
\usepackage{xcolor}
\usepackage{tikz}
\usetikzlibrary{matrix}
\usepackage{xargs}
\usepackage[all]{xy}
\usepackage{float}

\usepackage[margin=3cm]{geometry}

\usepackage{hyperref}
\hypersetup{
  hidelinks,
  colorlinks  = true,    
  urlcolor    = blue,    
  linkcolor   = blue,    
  citecolor   = blue,      
  pdfauthor   = { },%
  pdfsubject  = {Maniplexes},%
  pdftitle    = {},  %
}
\usepackage[capitalise,noabbrev, nameinlink,]{cleveref}

\theoremstyle{plain}
\newtheorem{thm}{Theorem}[section]
\newtheorem{theorem}[thm]{Theorem}

\newtheorem{lemma}[thm]{Lemma}

\newtheorem{proposition}[thm]{Proposition}

\newtheorem{corollary}[thm]{Corollary}
\newtheorem{definition}[thm]{Definition}
\newtheorem{remark}[thm]{Remark}

\crefname{proposition}{proposition}{propositions}

\newenvironment{red}{\relax\color{red}}{\hspace*{.5ex}\relax}
\newcommand{\bcp}{\begin{red}}
\newcommand{\ecp}{\end{red}}

\newenvironment{blue}{\relax\color{blue}}{\hspace*{.5ex}\relax}
\newcommand{\bgc}{\begin{blue}}
\newcommand{\egc}{\end{blue}}

\newenvironment{purple}{\relax\color{purple}}{\hspace*{.5ex}\relax}
\newcommand{\bpt}{\begin{purple}}
\newcommand{\ept}{\end{purple}}

\newcommand{\cC}{\mathcal{C}}
\newcommand{\cH}{\mathcal{H}}
\newcommand{\cK}{\mathcal{K}}

\newcommand{\cP}{\mathcal{P}}

\newcommand{\cT}{\mathcal{T}}

\newcommand{\G}{\Gamma}

\DeclareMathOperator{\aut}{Aut}
\DeclareMathOperator{\Ch}{Ch}

\usepackage[loadshadowlibrary,shadow, colorinlistoftodos,textsize=tiny,obeyFinal]{todonotes}

\setuptodonotes{fancyline, backgroundcolor=gray!10,bordercolor=gray}

\newcommandx{\inline}[2][1=]{\todo[inline, #1]{#2}}
\newcounter{hw}
\newcommandx{\homework}[2][1=]{\todo[inline, caption={Homework \thehw} #1]{\stepcounter{hw} #2}}

\newcommandx{\gabeLong}[2][1=Long todo, usedefault]{\todo[inline, color=blue!25, caption={\textbf{Gabe:} #1}]{\textbf{Gabe: #1}. #2}}
\newcommandx{\claudioLong}[2][1=Long todo, usedefault]{\todo[inline, color=red!25,caption={\textbf{Claudio:} #1}]{\textbf{Claudio: #1}. #2}}
\newcommandx{\philippeLong}[2][1=Long todo, usedefault]{\todo[inline,color=purple!25, caption={\textbf{Philippe:} #1}]{\textbf{Philippe: #1}. #2} }

\title{Flatness and Tightness in Hypertopes}
\author{Gabe Cunningham}
\address{Wentworth Institute of Technology \\
Boston, MA, USA 02115}

\author{Claudio Alexandre Piedade}
\address{Université Libre de Bruxelles \\
1050 Brussels, Belgium}
\author{Philippe Tranchida}
\address{Université Libre de Bruxelles \\
1050 Brussels, Belgium}

\email{cunninghamg1@wit.edu, claudio.piedade@ulb.be, tranchida.philippe@gmail.com}
\date{\today}

\begin{document}

\begin{abstract}
A \emph{hypertope} is a thin, residually connected incidence geometry. For each pair of types $\{i,j\}$, the residues of that type must be polygons, and if all residues of that type are the same type of polygon, then we can record this information in a diagram $D$. Given such a diagram, we prove a lower bound on the number of chambers of a hypertope with that diagram. A hypertope that achieves this bound is called \emph{tight}. We prove that a hypertope is tight if and only if its diagram $D$ is a forest and, for every pair of types $\{i,j\}$ that is not an edge of of $D$, every $i$-element is incident to every $j$-element. This generalizes a previously known result about abstract polytopes. Finally, we describe a combinatorial method and an algebraic method for constructing tight hypertopes.
\end{abstract}

\maketitle


\section{Introduction}

The theory of abstract polytopes is rich and by now well-established. It combines elements of geometry, combinatorics and group theory. Abstract polytopes are usually defined as posets satisfying a set of defining properties inspired by those of convex polytopes. From the point of view of incidence geometry, an abstract polytope is defined as a thin, residually connected geometry with a linear Buekenhout diagram \cite{buekenhout}.
The linearity of the diagram induces a partial order on the elements and allows for the use of many geometrical tools. Nonetheless, much of the underlying combinatorics and group theory remains valid regardless of the diagram. Hypertopes, first introduced in \cite{hypertopes}, are the objects obtained by dropping the linearity condition on the diagram of a polytope. They are thus thin, residually connected geometries.

A hypertope is called \emph{regular} when it has the highest degree of symmetry, meaning that its automorphism group acts transitively on its chambers. It turns out that regular hypertopes can always be constructed from smooth quotients of Coxeter groups. From this point of view, allowing nonlinear diagrams is a natural step. It brings into consideration geometries associated with arbitrary Coxeter groups, allowing for new behaviours and symmetries while preserving a considerable amount of the structure of polytopes. This broader setting also connects the theory to other well-known objects, such as hypermaps and Coxeter complexes. Moreover, nonlinear diagrams allow for other interactions between the different types of elements, and in particular for symmetry phenomena that cannot occur in the linear setting, such as trialities or higher order correlations \cite{HypertopesDiagramSymmetries,leemans2025groups}.

Even when focusing on polytopes, it can be beneficial to distinguish between properties that are intrinsically related to the linearity of the diagram, and thus to the natural partial order on the faces, and those that are more combinatorial or group-theoretical in nature. We therefore believe that studying hypertopes can help deepen our understanding of thin, residually connected geometries and thus, in particular, of abstract polytopes.

As with hypertopes, the study of abstract polytopes is primarily concentrated on those that are regular or otherwise highly symmetric. A natural way to remove this symmetry requirement while still insisting on some degree of homogeneity is to study \emph{equivelar} polytopes, which are precisely polytopes that have a Buekenhout diagram where every edge $\{i,j\}$ is labeled with a single number $p_{ij}$. In other words, these are the polytopes such that residues of each type $\{i,j\}$ ($i \neq j$) are $p_{ij}$-gons. Once we fix a diagram, a simple counting argument establishes a lower bound on the number of chambers (\cite[Proposition 3.3]{tight-polytopes}). A polytope that achieves this lower bound is called \emph{tight}. Several papers have explored the structure of tight polytopes, whether in the presence of high symmetry \cite{tight2,tight-chiral-polyhedra,tight3,tight-chiral-polytopes} or without it \cite{tight-polytopes}. There are two important general results. The first is that the diagram of a tight polytope cannot have two odd labels in a row \cite[Proposition 3.9]{tight-polytopes}. The second is that an $n$-polytope is tight if and only if it is $(i,i+2)$-flat for each $0 \leq i \leq n-3$ (that is, every $i$-element is incident to every $(i+2)$-element) \cite[Theorem 4.4]{tight-polytopes}. 

In this paper, we seek to extend the machinery from \cite{tight-polytopes} to hypertopes. First, we give a lower bound on the number of chambers of a hypertope with a given diagram $D$. A hypertope that meets this lower bound is called \emph{tight}. We investigate two questions:
\begin{enumerate}
    \item If a hypertope is tight, what restrictions does that put on its diagram $D$?
    \item Is there a relation between tightness and flatness, analogous to the situation with polytopes?
\end{enumerate}
This culminates in Theorem~\ref{tight-iff-flat}, where we find that a hypertope is tight if and only if its diagram $D$ is a forest and, for every non-edge $\{i,j\}$ of $D$, every $i$-element is incident to every $j$-element.

We also describe two constructions for building tight hypertopes. The first is combinatorial, adapting the techniques of \cite{flat-extensions}. Indeed, the construction more generally can be used to build ``small'' extensions of an existing hypertope, and the properties of this extension are collected in Theorem~\ref{flat-ext-thm}. The second construction is algebraic and uses the techniques described in \cite{WreathCosetGeom} to build tight regular hypertopes. 

We start by giving background on incidence systems and hypertopes in Section~\ref{sec:bkgd}. Then we generalize the notion of \emph{flatness} from polytopes to incidence systems in Section~\ref{sec:flat}. In Section~\ref{sec:tight}, we establish a lower bound on the size of a hypertope with a given diagram, and we derive several structural results on those hypertopes that meet the lower bound (\emph{tight} hypertopes). Finally, we describe constructions for flat and tight hypertopes in Sections~\ref{sec:flat-ext} and~\ref{sec:twisting}.

\section{Background}
\label{sec:bkgd}

\subsection{Incidence systems}
Incidence systems are a mathematical concept that encodes a great many of the geometrical and combinatorial structures that one may encounter in mathematics. We give here a brief account of the main definitions and terminology of incidence systems.
For a more detailed introduction to incidence geometry, we refer to~\cite{buekenhout}.

    A triple $\Gamma = (X,*,t)$ is called an \textit{incidence system} over $I$ if
    
\begin{itemize}
        \item $X$ is a set whose elements are called the \textit{elements} of $\Gamma$,
        \item $*$ is a symmetric and reflexive relation (called the \textit{incidence relation}) on $X$, and
        \item $t$ is a map from $X$ to $I$, called the \textit{type map} of $\Gamma$, such that distinct elements $x,y \in X$ with $x * y$ satisfy $t(x) \neq t(y)$. 
\end{itemize}

Elements of $t^{-1}(i)$ are called the elements of \emph{type $i$}.
The \textit{rank} of $\Gamma$, denoted as $\textnormal{rank}(\Gamma)$, is the cardinality of the type set $I$, which we will assume is finite for the remainder of the paper.
A \textit{flag} in an incidence system $\Gamma$ over $I$ is a set of pairwise incident elements. The type of a flag $F$ is $t(F)$, that is the set of types of the elements of $F$. The \emph{rank} of a flag $F$ is the cardinality of $t(F)$, while is \emph{corank} is $\textnormal{rank}(\Gamma)-|t(F)|$. A \textit{chamber} is a flag of type $I$. An incidence system $\Gamma$ is an \textit{incidence geometry} if all its maximal flags are chambers. 
We will denote the set of chambers of $\Gamma$ by $\Ch(\G)$. 

Let $F$ be a flag of $\Gamma$. An element $x\in X$ is {\em incident} to $F$ if $x*y$ for all $y\in F$. The \textit{residue} of $\Gamma$ with respect to $F$, denoted by $\Gamma_F$, is the incidence system formed by all the elements of $\Gamma$ incident to $F$ but not in $F$. The \textit{rank} of the residue $\Gamma_F$ is equal to rank$(\Gamma)$ - $|t(F)|$, the corank of $F$. The \emph{type} of a residue $\Gamma_F$ is the set $I\setminus t(F)$; equivalently, it is the cotype of $F$. Thus, a residue of type $J$ is the residue of a flag of cotype $J$. Said another way, the type of a residue is the set of types that remains.

The \textit{incidence graph} of $\Gamma$ is a graph with vertex set $X$ and where two elements $x$ and $y$ of different types are connected by an edge if and only if $x * y$.
An incidence geometry $\Gamma$ is \textit{connected} if its incidence graph is connected. It is \textit{residually connected} if all its residues of rank at least two are connected.
An incidence geometry is \textit{thin} if all its residues of flags of corank one contain exactly two elements.

A thin residually connected incidence geometry is called a \textit{hypertope}. Unless otherwise stated, we will assume that our hypertopes are finite.

An {\em automorphism} of an incidence geometry $\Gamma$ is a bijection $\phi$ of $X$ respecting the incidence relation $*$ and such that $t(\phi(x)) = t(x)$ for all $x \in X$. The group of all automorphisms of $\Gamma$ is denoted by $\aut(\Gamma)$.

The incidence geometry $\Gamma$ is said to be \textit{flag-transitive} if the action of $\aut(\Gamma)$ on the set of flags of each type 
is transitive. It is said to be \textit{incidence-transitive}
if, for each $i,j \in I$, the action of $\aut(\G)$ is transitive on pairs $(x, y)$ such that $t(x) = i$, $t(y) = j$, and $x * y$. A flag-transitive hypertope is also called \emph{regular}.

\subsection{Coset incidence systems}
Highly symmetric incidence geometries can be constructed from a group $G$ together with a set $(G_i)_{i \in I}$ of subgroups of $G$ as described in~\cite{Tits1957}. 

    A \emph{coset incidence system} $\Gamma(G,(G_i)_{i \in I})$, denoted also as $(G,(G_i)_{i \in I})$, is the incidence system over the type set $I$ where:
    \begin{enumerate}
        \item The elements of type $i \in I$ are left cosets of the form $g G_i $, $g \in G$.
        \item The incidence relation is given by non-empty intersection. More precisely, the element $g G_i$ is incident to the element $k G_j $ if and only if $g G_i \cap k G_j \neq \emptyset$.
    \end{enumerate}
We usually refer to the subgroups $G_i$ as the \emph{(standard) maximal parabolic subgroups}. For $J\subseteq I$, we define \emph{(standard) parabolic subgroups} as $G_J=\cap_{i\in J} G_i$, and define $G_\emptyset=G$.
When $J=I\setminus\{i\}$, we say that $G_J=G_{I\setminus\{i\}}$ are the \emph{minimal parabolic subgroups}, and we will denote them as $G^i$. Finally, the subgroup $G_I=\cap_{i\in I} G_i$ is called the \emph{Borel subgroup} of $\Gamma$. 

Any incidence system that contains a chamber and is incidence-transitive can be constructed as a coset incidence system.

\begin{lemma}\cite[Lemma 1.8.2]{buekenhout}
\label{inc-trans-coset-geo}
    Let $\Gamma$ be an incidence system such that $G =\aut(\Gamma)$ acts incidence-transitively on $\Gamma$. Let $\{x_i| i \in I\}$ be a chamber of $\Gamma$. Then, we have $\Gamma \cong (G,(G_i)_{i \in I})$ where $G_i$ is the stabilizer of $x_i$ in $G$.
\end{lemma}

\subsection{Chamber systems}\label{sec:chambersystems}

The concept of chamber system was invented by Tits in~\cite{tits1982}.
A \textit{chamber system} over $I$ is a pair $\cC=(C,\{\sim_i,i\in I\})$ consisting of a set $C$, whose members are called \textit{chambers}, and a collection of equivalence relations $\sim_i$ on $C$, indexed by $i \in I$.
Two chambers $c$ and $d$ are called $i-$\textit{adjacent} if $c \sim_i d$. For $i\in I$, each $\sim_i$-equivalence class is called an \textit{i-panel}.
The \textit{rank} of $\cC$ is $|I|$. The chamber system $\cC$ is called \textit{thin} if every $i-$panel is of size exactly $2$.

For $J \subset I$, we denote by $\sim_J$ the union of all $\sim_j$ with $j\in J$. A connected component of the graph of $(\cC,\sim_J)$ is called a \textit{J-cell} of $\cC$. For $i\in I$, the $(I$ \textbackslash $\{i\})$-cells are called \emph{$i$-elements} of $\cC$ (called \emph{$i$-objects} in \cite[Definition 3.2.2]{buekenhout}). In this article, we will only deal with thin chamber systems. In that case, the chamber system $\cC$ over $I = \{0,1,\cdots,n-1\}$ can be thought of as an $n$-valent graph with a proper $n$-coloring of its edges. The vertices of this graph correspond to the chambers of $\cC$ and two distinct chambers are connected by an $i$-colored edge if and only if they are $i$-adjacent. Whenever $\cC$ is thin, it guarantees that there is exactly one edge of each color at each vertex. In this context, a $J$-cell of $\cC$ is a connected component of the graph obtained by deleting all edges of color in $I -J$.

A chamber system $\cC$ over $I$ is \textit{residually connected} if, for every subset $J$ of $I$ and every system of $j$-elements $Z_j$, one for each $j \in J$, with the property that any two have a non-empty intersection, it follows that $\bigcap_{j \in J} Z_j$ is an $(I$ \textbackslash$J$)-cell.

It is always possible to construct 
an incidence system from a chamber system. Beware that in general, information can be lost in the process. If $\cC$ is a chamber system, the \textit{incidence system of $\cC$}, denoted by $\Gamma(\cC)$, is the incidence system over $I$ determined as follows. Its $i$-elements, for $i\in I$ are the pairs $(x,i)$ with $x$ an $i$-element of $\cC$; two elements $(x,k)$ and $(y,l)$ of $\Gamma(\cC)$ are incident if and only if $x \cap y \neq \emptyset$ in $\cC$, i.e., $x$ and $y$ have a chamber in common. (The pairs contain the type $i$ because it may happen that for some $i \neq j$, an $i$-element and a $j$-element consist of the same chambers, and we need to distinguish the corresponding elements in the incidence system.) 

The next theorem is a simplified version of Theorem 3.4.6 in \cite{buekenhout}. Please refer to the original textbook for the more precise statement, which requires a bit more background that what is given in this section.

\begin{thm}\label{thm:Correspondence} 
Let $\mathcal{G}(I)$ be the collection of all residually connected geometries over $I$ and let $\cC(I)$ be the collection of all residually connected chamber systems over $I$.
    There is a bijective homomorphism preserving structures between the collections $\mathcal{G}(I)$ and $\cC(I)$. Thus, each residually connected geometry over $I$ corresponds to a unique residually connected chamber system over $I$ (up to isomorphism), with the same automorphism group, and vice versa.
\end{thm}

Therefore, when considering only residually connected incidence geometries and residually connected chamber system, one can freely go back and forth between the two view points.

\subsection{Diagrams}

Francis Buekenhout introduced a diagram associated to incidence geometries ~\cite{buekenhout}. The main idea is that local information, given by rank two residues, often imposes global conditions on the whole geometry. The definition of diagram holds for any incidence geometry. Here, we will give the definition in the context of hypertopes only, for simplicity.

Let $\cH$ be a hypertope, that is, a residually connected and thin incidence geometry. We have seen in Section~\ref{sec:chambersystems} that we can also consider $\cH$ as a chamber system. Let $I$ be the type set of $\cH$. For any $i \neq j \in I$, since $\cH$ is thin, any rank two residue of type $\{i,j\}$ of $\cH$ must be the geometry of an $n$-gon, for some integer $n \geq 2$. From the point of view of thin chamber systems, the graph obtained from the chamber graph of $\cC(\cH)$ by deleting all edges of color not equal to $i$ or $j$ must be a cyclic graph of length $2 n$ with edges of alternating colors. Of course, the integer $n$ usually depends on the rank two residue we chose, and not only on the types $i$ and $j$. The integer $n$ is called the \textit{parameter} of the residue.

In this paper, a \textit{diagram} is a labeled graph $D$ whose vertex set is $I$ and whose edges are labeled by integers strictly greater than two. By convention, a missing edge in $D$ should be understood as an edge with label $2$. We say that $\cH$ \textit{belongs to} $D$, or that $\cH$ \textit{has diagram} $D$, if for every $i\neq j \in I$, all rank two residues of type $\{i,j\}$ of $\cH$ have the same parameter $p_{ij}$ (that is, the residues are $p_{ij}$-gons), 
and that this parameter is exactly the label of the edge $\{i,j\}$ of $D$.\footnote{In \cite[Definition 2.3.7]{buekenhout}, the authors allow diagrams where each edge is labeled by a whole class of parameters. In our usage, each edge has only a single label. Thus, when we say that a hypertope \emph{belongs to a diagram}, we are asserting that there is a diagram in our sense (with a single label per edge), and not in the more general sense of \cite{buekenhout}.}

Note that a hypertope belonging to a linear diagram is an abstract polytope (see \cite{arp}). In fact, it is an \emph{equivelar} polytope, which is to say that it has a \emph{Schl\"{a}fli symbol}. Thus, we may think of hypertopes that belong to a diagram $D$ as being a generalization of equivelar abstract polytopes.

Diagrams for more general incidence geometries are more detailed, because the rank two residues are much more varied. In this paper, our results that involve diagrams of incidence geometries only need information about the underlying graph, so will not worry about these additional parameters.

The diagram of a hypertope provides an easy lower bound on the number of $i$-elements for each $i$:

\begin{proposition}
\label{num-i-elements}
In a hypertope with diagram $D$, the number of $i$-elements is at least $\max_{j \neq i} p_{ij}$.
\end{proposition}

\begin{proof}
In order for the residue of type $\{i,j\}$ to be a $p_{ij}$-gon, there must be at least $p_{ij}$ $i$-elements.    
\end{proof}

In other words, if we look at node $i$ in the diagram, then there are at least as many $i$-elements as the largest label among all edges incident to $i$. Recall that if there is no edge between $i$ and $j$ in the diagram, we consider that the same as there being an edge of label $2$, so that $p_{ij} = 2$ in this case.


Note that if $\cH$ is regular, it is incidence-transitive, and thus can be constructed as a coset geometry $(G,(G_i)_{i \in I})$ (see Lemma~\ref{inc-trans-coset-geo}). The thinness of $\cH$, together with residual-connectedness, forces the minimal parabolic subgroups $G^i = \cap_{j \neq i}G_j$ to be cyclic groups of order two. Hence $G^i = \langle \rho_i \rangle$ for some involution $\rho_i \in G = \aut(\cH)$. Then, when $\cH$ is regular, the integer $p_{ij}$ can also be understood as the order of $\rho_i \rho_j$ in $G$. Therefore, in the context of highly symmetric hypertopes, the Buekenhout diagram of $\cH$ is equivalent to the Coxeter diagram of $G = \aut(\cH)$.

For some diagrams, the incidence relation gains some degree of transitivity (not to be confused with incidence-transitivity, which is about the action of $\aut(\G)$ on incident pairs). Though we are mainly interested in hypertopes, the following result applies to any incidence geometry with a diagram.

\begin{lemma}
\label{trans-incidence}
Let $\G$ be a residually connected incidence geometry with diagram $D$. Suppose that $x$ is an $i$-element, $y$ is a $j$-element, and $z$ is a $k$-element such that $x * y$ and $y * z$. If $i$ and $k$ are in separate connected components of $D \setminus j$, then $x * z$.
\end{lemma}

\begin{proof}
Both $x$ and $z$ are in the residue $\G_y$, and the diagram of $\G_y$ is a subgraph of $D \setminus j$, and so has $i$ and $k$ in separate components. It follows that every $i$-element of $\G_y$ is incident to every $k$-element by \cite[Theorem 2.1.6]{buekenhout}, and so $x * z$.
\end{proof}

\subsection{Direct sums}

There is a natural notion of the \emph{direct sum} of incidence geometries \cite[Definition 2.1.8]{buekenhout}. Specifically, suppose $\{I_j\}_{j \in J}$ is a collection of pairwise disjoint sets, and consider a collection of incidence geometries $\G_j = (X_j, *_j, t_j)$ over $I_j$, where the sets $X_j$ are pairwise disjoint. Then the direct sum of these is an incidence geometry $\G = (X, *, t)$ over $\bigcup_j I_j$, where $X = \bigcup_j X_j$ and the type function $t$ is defined by $t \mid_{X_j} = t_j$. Incidence is defined by requiring $* \mid_{X_j} = *_j$, and then setting $x * y$ for all $x \in X_i, y \in X_k$ for $i \neq k$. That is, incidence between two elements in the same $X_j$ remains unchanged, and every pair of elements from different sets $X_i$ and $X_k$ is incident. 

For $\G$ as defined above, if the individual summands $\G_j$ have diagram $D_j$, then $\G$ has diagram $\bigsqcup_j D_j$. The converse is true for hypertopes:

\begin{proposition}{\cite[Corollary 2.1.10]{buekenhout}}
\label{dir-sum}
If $\cH$ is a hypertope with diagram $D$, and $D$ has connected components $D_1, \ldots, D_m$, then $\cH$ is isomorphic to the direct sum of $\cH_1, \ldots, \cH_m$, where $\cH_j$ is the hypertope corresponding to $D_j$ (that is, the induced geometry with all elements of the types in $D_j$).
\end{proposition}

\section{Flatness}
\label{sec:flat}

    In the theory of abstract polytopes, a polytope is said to be $(i,j)$-flat if every $i$-face ($i$-element) is incident to every $j$-face (see \cite[Section 4E]{arp}). It makes sense to extend this definition to arbitrary incidence systems (though we will primarily be concerned with hypertopes):

    \begin{definition}
    For types $i \neq j$, an incidence system $\Gamma$ is \emph{$(i,j)$-flat} if every $i$-element is incident to every $j$-element. If an incidence system is $(i,j)$-flat for some choice of $i$ and $j$, then we also say simply that it is \emph{flat}.
    \end{definition}

    Equivalently, an incidence system is $(i,j)$-flat if the subgraph of the incidence graph that is induced by all nodes of types $i$ and $j$ is a complete bipartite graph. 

    \begin{proposition}
    \label{flat-hereditary}
    Suppose that $\G$ is an incidence system that is $(i,j)$-flat. Then the residues of $\G$ that contain types $i$ and $j$ are $(i,j)$-flat. Furthermore, if $\G$ is an incidence geometry with diagram $D$, then $i$ and $j$ are not adjacent in $D$. 
    \end{proposition}

    \begin{proof}
    The first part is clear: every $i$-element is incident to every $j$-element, and this incidence persists in the residue. When $\G$ is an incidence geometry, being $(i,j)$-flat means that the residues of type $\{i,j\}$ are \emph{generalized digons} (\cite[Definition 2.1.1]{buekenhout}). Thus, $i$ and $j$ are not adjacent in the diagram of $\G$.
    \end{proof}

    Next, we prove an analogue of \cite[Lemma 4E3]{arp} for hypertopes, which is a partial converse to Proposition \ref{flat-hereditary}.

    \begin{proposition}
    \label{inherit-flat}
    Let $\cH$ be a residually connected incidence geometry with diagram $D$. Consider types $i$, $j$, and $k$ such that $i$ and $k$ are in separate connected components of $D \setminus j$. If the residue of every $k$-element of $\cH$ is $(i,j)$-flat, then $\cH$ is $(i,j)$-flat.
    \end{proposition}

    \begin{proof}
    Let $x$ be an $i$-element and $y$ be a $j$-element of $\cH$. Since $\cH$ is residually connected, there is a sequence
    \[ x, z_1, y_1, z_2, y_2, \ldots, z_m, y_m = y\]
    such that each $y_a$ is a $j$-element, each $z_a$ is a $k$-element, and every pair of adjacent elements is incident (see \cite[Lemma 1.6.3]{buekenhout}). Now, $x$ is an $i$-element that is incident to $z_1$ and $y_1$ is a $j$-element that is incident to $z_1$ so that $x$ and $y_1$ are in the residue $\cH_{z_1}$. Since the residue at every $k$-element is $(i,j)$-flat, it follows that $x$ is incident to $y_1$. Furthermore, $y_1$ is incident to $z_2$. Then by Lemma~\ref{trans-incidence}, it follows that $x$ is incident to $z_2$. Then again by flatness of every $k$-residue, we get that $x$ is incident to $y_2$, and continuing in this way we establish that $x$ is incident to $y$, and so $\cH$ is $(i,j)$-flat.
    \end{proof}

	\begin{proposition}
	\label{forced-flat}
	Let $\cH$ be an $n$-hypertope with diagram $D$. Suppose that for some $i$ and $j$, the number of $i$-elements is exactly
	$p_{ij}$. Then:
	\begin{enumerate}
	\item $\cH$ is $(i,k)$-flat for each $k \neq i,j$.
	\item $p_{ik} = 2$ for each $k \neq i,j$. In particular, in the diagram $D(\cH)$, $i$ is either an isolated node or a leaf node connected only to type $j$. 
	\end{enumerate}
	\end{proposition}
	
	\begin{proof}
	Suppose the number of $i$-elements is $p_{ij}$ and let $x$ be a $k$-element of $\cH$. For any flag $F$ of type $I \setminus \{i,j\}$ containing $x$, the residue $\cH_F$ is
	a $p_{ij}$-gon. Hence, all $i$-elements appear in this residue and therefore, we get that $x$ is incident to all $i$-elements of $\cH$. Therefore, every $k$-element is incident to every $i$-element, so $\cH$ is $(i,k)$-flat.
	
	The second part follows from the first and Proposition~\ref{flat-hereditary}.
	\end{proof}

    \begin{corollary}
    \label{flat-elt-count}
    Let $\cH$ be an $n$-hypertope with diagram $D$. The number of $i$-elements of $\cH$ is $p_{ij}$ if and only if $\cH$ is $(i,k)$-flat for each $k \neq i,j$. 
    \end{corollary}

    \begin{proof}
    Proposition~\ref{forced-flat} shows one direction. For the other, suppose $\cH$ is $(i,k)$-flat for each $k \neq i,j$, and consider a residue of type $\{i,j\}$. Since in $\cH$ every $i$-element is incident to every $k$-element for every $k \in I \setminus \{i,j\}$, it follows that every $i$-element is in this residue. On the other hand, the residue is a $p_{ij}$-gon, so there are $p_{ij}$ exactly $i$-elements.
    \end{proof}

\section{Tightness in hypertopes}
\label{sec:tight}

    Our next goal will be to establish a lower bound on the size of hypertopes with a given diagram, and to further characterize those hypertopes that meet the lower bound. Let us start by reviewing the situation for polytopes (hypertopes with a linear diagram).

    An equivelar polytope $\cP$ of type $\{p_1, \cdots, p_{n-1}\}$ has at least $2p_1 \cdots p_{n-1}$ chambers \cite[Proposition 3.3]{tight-polytopes}, and when $\cP$ has exactly that many chambers, it is called \emph{tight}. In rank 3 there is an obvious necessary condition that follows from the fact that the number of chambers must be a multiple of 4, and this condition turns out to be sufficient: 

    \begin{theorem}[{\cite[Theorem 5.1]{tight-polytopes}}]
    \label{tight-polyhedra}
    There is a tight polyhedron ($3$-polytope) of type $\{p, q\}$ if and only if $p$ and $q$ are not both odd.
    \end{theorem}
    
    If we are looking for tight \emph{regular} (flag-transitive) polyhedra, then the situation is more nuanced. The result below follows from \cite[Theorem 3.4]{tight2} and \cite[Theorem 5.10]{tight3}.

    \begin{theorem}
    \label{tight-regular-polyhedra}
    There is a tight regular polyhedron of type $\{p, q\}$ if and only if one of the following is true:
    \begin{enumerate}
        \item $p$ and $q$ are even.
        \item $p$ is odd and $q$ is an even divisor of $2p$ or $q$ is odd and $p$ is an even divisor of $2q$.
        \item $p$ is an odd multiple of $3$ and $q = 4$ or $q$ is an odd multiple of $3$ and $p = 4$.
    \end{enumerate}
    \end{theorem}



    The lower bound for the size of equivelar polytopes follows from a simple chamber-counting argument, so let us try the same idea for hypertopes. To begin, we make the following definition: if $T$ is a tree with edge labels, then the \emph{weight} of $T$, denoted $w(T)$, is the product of the edge labels. The weight of a single isolated vertex is defined to be 1. The weight of a forest $\cT = T_1 \sqcup \cdots \sqcup T_m$ is defined to be $w(\cT) = 2^{m-1} \prod_{i=1}^m w(T_i)$. Thought of another way, given a forest $\cT$ consisting of $m$ trees, we can turn it into a tree by adding $m-1$ edges with label $2$, and then $w(\cT)$ is the weight of this new tree.

	\begin{proposition}
    \label{flag-bound}
    If $\cH$ is a hypertope with connected diagram $D$, and $T$ is a maximum-weight spanning tree of $D$, then $\cH$ has at least $2 w(T)$ chambers.
	\end{proposition}
	
	\begin{proof}
    It suffices to prove the bound for \emph{all} spanning trees of $D$. When $n = 1$, the diagram $D$ is a single node, which has weight 1, and $\cH$ must be the geometry of an edge, with two chambers. When $n = 2$, $\cH$ must be a polygon and the result is also true. Now suppose $n \geq 3$ and that the result is true for $(n-1)$-hypertopes. Let $T$ be a spanning tree of $D$, and let $j$ be one of the leaf nodes. Suppose an edge of weight $p$ connects $j$ to $T \setminus j$. 
	Then we may specify a chamber of $\cH$ by picking an element $x$ of type $j$, and then picking a chamber in the residue of $x$. There are at least $p$ choices of element $x$ (by Proposition~\ref{num-i-elements}), and since $T \setminus j$ is a spanning tree of $D\setminus j$,
	the residue $\cH_x$ has at least $2 w(T\setminus j)$ chambers, by inductive hypothesis. The result follows.
	\end{proof}

    We will be discussing \emph{spanning forests} of a diagram $D$ in this section, by which we will always mean a collection of spanning trees, one per connected component of $D$. (For example, we are not considering taking a spanning tree of a connected graph and then just removing one edge.)  

    \begin{corollary}
    \label{flag-bound-forest}
    If $\cT = T_1 \sqcup \cdots \sqcup T_m$ is a spanning forest of the diagram $D$ of $\cH$, then $|\Ch(\cH)| \geq 2w(\cT)$.
    \end{corollary}

    \begin{proof}
    Since $\cT$ is a spanning forest with $m$ trees, $\cH$ has diagram $D = D_1 \sqcup \cdots \sqcup D_m$, where $T_i$ is a spanning tree of $D_i$ for each $i$. Then, since $\cH$ is residually connected, $\cH$ is the direct sum of the geometries corresponding to each $D_i$ (by Proposition~\ref{dir-sum}), and the number of chambers is at least the product of the number of chambers from each geometry. So the number of chambers is at least
    \[ \prod_{i=1}^m 2w(T_i) = 2 w(\cT). \]
    \end{proof}

    Notice that if $\cH$ has a linear diagram, then this is precisely the same bound as what was already known for polytopes. 

    With this new size bound for hypertopes, we make the following definition:

	\begin{definition}
	If $\cH$ is a hypertope with diagram $D$, we say that $\cH$ is \emph{tight} if $|\Ch(\cH)| = 2 w(\cT)$ for a maximum-weight spanning forest $\cT$ of $D$.
	\end{definition}

    \begin{proposition}
    \label{tight-discon}
    Suppose $\cH$ has diagram $D_1 \sqcup \cdots \sqcup D_m$. Then $\cH$ is the direct sum of hypertopes $\cH_1, \ldots, \cH_m$, and $\cH$ is tight if and only if each $\cH_i$ is tight. 
    \end{proposition}

    \begin{proof}
    That $\cH$ is a direct sum follows from Proposition~\ref{dir-sum}. Furthermore, since every element of $\cH_i$ is incident to every element of $\cH_j$ whenever $i \neq j$, it follows that
    \[ |\Ch(\cH)| = |\Ch(\cH_1)| \cdots |\Ch(\cH_m)|. \]
    Then $\cH$ is tight if and only if $|\Ch(\cH)| = 2w(\cT) = 2^m \prod_{i=1}^m w(T_i)$, where $\cT = T_1 \sqcup \cdots \sqcup T_m$ is a max weight spanning forest of $\cH$. Now, Proposition~\ref{flag-bound} says that
    \[ |\Ch(\cH_1)| \cdots |\Ch(\cH_m)| \geq 2^m \prod_{i=1}^m w(T_i), \]
    and we get equality only when $|\Ch(\cH_i)| = 2w(T_i)$ for each $i$, and the result follows. 
    \end{proof}
    
	\begin{proposition} \label{tight-residues}
	Suppose that $\cH$ is a tight $n$-hypertope with $n \geq 3$, and let $i$ be a leaf of a maximum-weight spanning forest $\cT$.
	\begin{enumerate}
	\item[(a)] If $i$ is isolated, then there are two $i$-elements, and otherwise the number of $i$-elements is equal to the weight of the edge at $i$.
	\item[(b)] The residue of $\cH$ at every $i$-element is tight.
	\item[(c)] For every type $k \neq i$ not adjacent to $i$ in $\cT$, the hypertope $\cH$ is $(i,k)$-flat.
	\end{enumerate}
	\end{proposition}
	
	\begin{proof}
    Let $n_i$ be the number of $i$-elements. Let $p$ be the label of the unique edge at $i$, or $2$ if there is no label.
    We can pick a chamber of $\cH$ by first picking an $i$-element $x$ and then picking a chamber in the residue $\cH_x$. By Proposition~\ref{flag-bound-forest}, the residue $\cH_x$ has at least $2 w(\cT \setminus i)$ chambers, so
    \[ |\Ch(\cH)| \geq n_i 2w(\cT \setminus i). \]
    On the other hand, by assumption we have
    \[ |\Ch(\cH)| = 2w(\cT) = 2 w(\cT\setminus i) p, \] 
    and so $n_i \leq p$. Since we also have $n_i \geq p$ by Proposition~\ref{num-i-elements}, it follows that $n_i = p$.
    Then it also follows that $|\Ch(\cH_x)| = 2w(\cT \setminus i)$, which means that $\cH_x$ is tight. 
    The last part follows from the first part and from Proposition~\ref{forced-flat}. 
	\end{proof}

    Now we come to our first big structural result about tight hypertopes.

	\begin{theorem}
	\label{tree-diagram}
    If $\cH$ is a tight hypertope with connected diagram $D$, then $D$ is a tree.
	\end{theorem}
	
	\begin{proof}
	We prove the result by induction on $n$, the rank of the hypertope. It is trivially true for $n = 2$.
	Suppose that $n \geq 3$ and that the result is true for $(n-1)$-hypertopes. Let $T$ be a maximum-weight spanning tree for the diagram $D$. Let $i$ be one of the leaf nodes of $T$, and let $j$ be the neighbor
	of $i$ in $T$. Then by Proposition~\ref{tight-residues}(c), the hypertope $\cH$ is $(i, k)$-flat for each $k \neq i,j$,
	and by Proposition~\ref{flat-hereditary}, that implies that $p_{ik} = 2$ for $k \neq i,j$. In other words, the only
	proper edge in $D$ that involves $i$ is the edge between $i$ and $j$. Now, by Proposition~\ref{tight-residues}(b),
	any residue $\cH_x$ with $x$ an $i$-element is tight, and so by inductive hypothesis, the diagram $D \setminus i$ is a tree. Since $D$ is obtained from
	$D \setminus i$ by adding a leaf node that is attached to a tree, no new cycles are created, and so $D$
	is itself a tree.
	\end{proof}

	\begin{corollary}
    \label{tight-3-hypertopes}
	All tight 3-hypertopes are polytopes (polyhedra).
	\end{corollary}

    \begin{proof}
    A tight 3-hypertope either has a disconnected diagram, or a tree diagram which must then be linear, and a hypertope with such a diagram is precisely an abstract polytope.
    \end{proof}

    Now we can extend \cref{tight-residues} to show that \emph{every} residue of a tight hypertope is tight.

    \begin{corollary}
    \label{all-residues-tight}
    If $\cH$ is a tight $n$-hypertope with $n \geq 3$, then every residue of $\cH$ is tight. 
    \end{corollary}

    \begin{proof}
    It suffices to prove that the residue at any element $x$ of type $i$ is tight. By Proposition~\ref{tight-discon} and Theorem~\ref{tree-diagram}, we may assume that the diagram $D$ of $\cH$ is a forest. Let $D_1 \sqcup \cdots \sqcup D_m$ be the connected components of the diagram $D \setminus i$. Since $\cH$ is residually connected, $\cH_x$ is a direct sum of hypertopes $\cH_1, \ldots, \cH_m$ with these diagrams. For each $r$, we can obtain $D_r$ by successively deleting leaves from $D$, and so Proposition~\ref{tight-residues}(b) proves that each $\cH_r$ is tight. It follows that $\cH_x$ itself is tight.
    \end{proof}

    All together, Theorem~\ref{tree-diagram}, Corollary~\ref{tight-3-hypertopes}, and Corollary~\ref{all-residues-tight} tell us that tight hypertopes are made by suitably gluing together tight polyhedra such that the diagram is a tree. Thus, restrictions on tight polyhedra yield restrictions on tight hypertopes of all ranks.

    \begin{proposition}
    If $\cH$ is a tight hypertope with connected diagram $D$, then $D$ is a tree with no adjacent odd labels.
    \end{proposition}

    \begin{proof}
    Two adjacent edges in $D$ determine a residue that must be a tight polyhedron, and the result follows from Theorem~\ref{tight-polyhedra}.
    \end{proof}

    \begin{proposition}
    If $\cH$ is a tight regular hypertope with connected diagram $D$, then each edge with an odd label $p$ is only adjacent to edges with labels that are even divisors of $2p$, unless $p$ is a multiple of $3$ and is adjacent only to edges with labels $6$ and $4$.
    \end{proposition}

    \begin{proof}
    This follows from \cref{tight-regular-polyhedra} considering that every rank 3 residue of a tight regular hypertope must be a tight regular polyhedron.
    \end{proof}

    Having worked out several necessary conditions for tight hypertopes, we will describe some constructions for them in the next section. First, though, let us work out more of the combinatorial information of tight hypertopes. First, we determine the number of elements of each type. 
    
    \begin{proposition}
    \label{element-counts}
    Suppose $\cH$ is a tight $n$-hypertope with diagram $T$. Then the number of $i$-elements is 
    \[ 2^{1-\deg_T(i)} \prod_{j: p_{ij}>2} p_{ij}. \]
    \end{proposition}

    \begin{proof}
    If $T$ is disconnected, so that $\cH$ is a direct sum (see Proposition~\ref{dir-sum}), then the number of $i$-elements of $\cH$ is the same as the number of $i$-elements of the hypertope $\cH'$ that corresponds to the connected component of $i$ in $T$. Thus, we may assume that $T$ is connected. Then since $\cH$ is tight, $T$ is a tree, and $|\Ch(\cH)| = 2w(T)$. On the other hand, we can specify a chamber by first picking an $i$-element $x$, and then picking a chamber in the residue of $x$. Thus
    \[ 2w(T) = n_i |\Ch(\cH_x)|, \]
    where $n_i$ is the number of $i$-elements. 
    Now, $\cH_x$ is tight by Corollary~\ref{all-residues-tight}. The diagram of $\cH_x$ consists of $d := \deg_T(i)$ connected components, trees $T_1 \sqcup \cdots \sqcup T_d$, so we have
    \[ |\Ch(\cH_x)| = 2w(T \setminus i) = 2^{d} \prod_{i=1}^d w(T_i). \]
    Now, 
    \[ \prod_{i=1}^d w(T_i) = \frac{w(T)}{\prod_{j : p_{ij}>2} p_{ij}}, \]
    and so
    \[ n_i = \frac{2w(T)}{|\Ch(\cH_x)|} = \frac{2w(T)}{2^{d} \prod_{i=1}^d w(T_i)} = 2^{1-d} \prod_{j: p_{ij}>2} p_{ij}. \]
    \end{proof}

    Now we can relate tightness to flatness in a way that extends \cite[Theorem 4.4]{tight-polytopes}. 

    \begin{theorem}
    \label{tight-iff-flat}
    A hypertope $\cH$ is tight if and only if it has a forest diagram $\cT$ and, for every non-edge $\{i,j\}$ of $\cT$, $\cH$ is $(i,j)$-flat.   
    \end{theorem}

    \begin{proof}
    We prove the result by induction on the rank $n$. The result is trivial for $n = 2$, since every polygon is tight. Suppose $n \geq 3$ and that the result is true for $(n-1)$-hypertopes. 
    
    First, suppose $\cH$ is tight. Proposition~\ref{tight-discon} and Theorem~\ref{tree-diagram} establish that $\cH$ has a forest diagram $\cT$. Now consider types $i$ and $j$ not joined by an edge in $\cT$.
    If $i$ and $j$ are in separate connected components of $\cT$, by \cite[Theorem 2.1.6]{buekenhout}, 
    every $i$-element is incident to every $j$-element, so $\cH$ is $(i,j)$-flat. Otherwise, we may suppose that $i$ and $j$ are in the same connected component.
    
    If either $i$ or $j$ is a leaf node, then Proposition~\ref{tight-residues}(c) says that $\cH$ is $(i,j)$-flat. Otherwise, there is some node $k$ that is adjacent to $j$ and not on the unique path between $i$ and $j$. Now, for every $k$-element $x$, the residue $\cH_x$ is a tight hypertope (by Corollary~\ref{all-residues-tight}). Since the diagram of $\cH_x$ is $\cT \setminus k$, and there is no edge $(i,j)$ in $\cT \setminus k$, it follows that $\cH_x$ is $(i,j)$-flat by inductive hypothesis.
    Note that removing node $j$ from $\cT$ would disconnect $i$ from $k$; otherwise there would be a cycle. Thus Proposition~\ref{inherit-flat} implies that $\cH$ itself is $(i,j)$-flat.

    Conversely, suppose that $\cH$ has a forest diagram $\cT$ and is $(i,j)$-flat whenever $(i,j)$ is not an edge of $\cT$. Consider a leaf node $i$ of $\cT$, connected by an edge of label $p$ to some type $j$. By assumption, $\cH$ is $(i,k)$-flat for each $k \in I \setminus \{i,j\}$, and so by Corollary~\ref{flat-elt-count}, the number of $i$-elements is $p$. Now, note that the residue at any given $i$-element $x$ also satisfies our hypotheses, and so is tight by inductive hypothesis. It follows that 
    \[ |\Ch(\cH)| = p |\Ch(\cH_x)| = 2p w(\cT \setminus i) = 2w(\cT), \]
    with the last equality following from the fact that the only edge weight missing from $\cT \setminus i$ compared to $\cT$ is the edge of weight $p$ at type $i$. Thus, $\cH$ is tight.
    \end{proof}

\section{Flat extensions of hypertopes}
\label{sec:flat-ext}

	Let $\cH$ be an $n$-hypertope. We say that the $(n+1)$-hypertope $\cK$ is an \emph{extension}
	of $\cH$ if there is a type $i$ such that all residues of $\cK$ of cotype $i$ are isomorphic to $\cH$.
	Suppose we start with $\cH$ and want to build an extension $\cK$. How can we accomplish this?
	
	One way that always works is to add a new type $i$ to the types $I$ of $\cH$, make exactly two
	elements of that type, and then make those two elements incident to every other element of $\cH$. Indeed, this is merely the direct sum of $\cH$ with a 1-hypertope. It is straightforward
	to show that this yields a hypertope. We will call hypertopes obtained this way \emph{trivial extensions}.
	Note that, due to thinness, we cannot try to do the same thing with three or more elements (that is, we cannot
	add three new elements and make them incident to every element of $\cH$).
	
	What other ``small'' extensions can we try? Suppose we want $\cK$ to have 4 $i$-elements $x$ such that $\cK_x$ is isomorphic
	to $\cH$. While we cannot make these elements incident to all other elements of $\cH$, could we pick a single
	type $j$ of $\cH$ and make each $i$-element incident to all elements of $\cH$ other than the $j$-elements? In other words,
	we want $\cK$ to be $(i,k)$-flat for $k \neq i,j$. By \cref{flat-hereditary}, this requires that $p_{ik} = 2$, which is
	to say that type $i$ is not connected to type $k$ in the diagram of $\cK$. So the diagram of $\cK$ consists of the
	diagram of $\cH$, with a single node added and connected to a single type $j$ of $\cH$ with
	an edge labeled $p_{i,j} \geq 3$. We will call an extension of this type a \emph{flat extension of $\cH$}. Without loss of generality, we will usually assume that $\cH$ has types $\{0, \ldots, n-1\}$ and that $j = n-1$. (Though let us emphasize that the linear order of types has no bearing here.)

    In \emph{Flat extensions of abstract polytopes}, the first author considered this problem in the context of polytopes \cite{flat-extensions}. The idea was to take several copies of the flag graph (chamber system) of the $n$-polytope $\cH$ and then add a new perfect matching corresponding to the new type. The new polytope had facets (the residues of the new elements) isomorphic to $\cH$, and it was $(n-2,n)$-flat. The same construction works here with very little modification; let us describe it.

    Fix an $n$-hypertope $\cH$ with chamber system $\cC$. For $i \in \{0, \ldots, n-1\}$, let $r_i$ be the permutation of the chambers of $\cC$ induced by edges of color $i$. If $C$ is a chamber of $\cH$ and $w = \langle r_0, \ldots, r_{n-1} \rangle$, we will use $C^w$ to denote the image of $C$ under $w$. We start with $k$ copies of $\cC$, with $k \geq 2$ an even integer, so that our chambers can be written as $(C, j)$ with $j \in \mathbb Z_k$ being the \emph{layer} of the chamber. We want to add edges of color $n$ such that the resulting chamber system is a hypertope and such that the diagram is the diagram of $\cH$ but with an extra leaf node of type $n$ connected to type $n-1$ with a branch of label $k$. Then the induced permutation $r_n$ must commute with $r_i$ for $i < n-1$ since the $\{i, n\}$-residue must be a digon for each $i < n-1$. This implies that once we determine where $r_n$ sends a chamber $(C, j)$, we know where it sends $(C^w, j)$ for every $w \in \langle r_0, \ldots, r_{n-2} \rangle$. 

    Let us now follow the idea of the construction in \cite[Section 3.1]{flat-extensions}. Suppose that $\cH$ is \emph{$(n-1)$-element bipartite}, meaning that it is possible to color the chambers of $\cC$ with colors red and blue such that adjacent chambers are different colors if and only if they are $(n-1)$-adjacent. 
    Color each copy of $\cC$ according to such a (single) coloring. Then we will match $(C,2m)$ with $(C,2m+1)$ if $C$ is red, and we will match $(C,2m)$ with $(C,2m-1)$ otherwise. We will denote this new graph by $\cC|k$, and the resulting incidence system by $\cH|k$. Clearly $\cC|k$ is a chamber system with $k$ copies of $\cC$ as its 
    $n$-elements. 
    We need to show that the result is a hypertope with the desired diagram.

    \begin{proposition}
    If $k \geq 4$ is even and $\cC$ is the chamber system of a hypertope, then so is $\cC|k$. 
    \end{proposition}

    \begin{proof}
    Clearly $\cC|k$ is thin, so we need only to show that it is residually connected. To do so, we make use of \cite[Lemma 3.4.9(ii)]{buekenhout}.  Let $Z_J$ ($Z_K, Z_L$) be a $J$-cell ($K$-cell, $L$-cell) such that they have pairwise nonempty intersection; we need to show that $Z_J \cap Z_K \cap Z_L$ is a $(J \cap K \cap L)$-cell. Set $\bar{J} = J \setminus \{n\}$ and similarly define $\bar{K}$ and $\bar{L}$. Note that since $r_n$ only affects the second coordinate of a chamber $(C, a)$, a $J$-cell of $\cC|k$ induces a $\bar{J}$-cell of $\cC$. Indeed, every $J$-cell $Z_J$ is contained in a set of the form $Z_{\bar{J}} \times \mathbb Z_k$, where $Z_{\bar{J}}$ is a $\bar{J}$-cell.

    First, suppose $n \not \in J \cap K \cap L$, and suppose without loss of generality that $n \not \in J$. Then $Z_J$ is a subset of $\cC \times \{a\}$ for some $a$, and since $Z_J$ has nonempty intersection with $Z_K$ and $Z_L$, then these cells also intersect $\cC \times \{a\}$. It follows that $Z_J \cap Z_K \cap Z_L = (Z_{\bar{J}} \cap Z_{\bar{K}} \cap Z_{\bar{L}}) \times \{a\}$. Since $\cC$ is a hypertope, $Z_{\bar{J}} \cap Z_{\bar{K}} \cap Z_{\bar{L}}$ is a $(\bar{J} \cap \bar{K} \cap \bar{L})$-cell of $\cC$, which is a $(J \cap K \cap L)$-cell of $\cC|k$ since $n \not \in J \cap K \cap L$.
    
    Next, suppose $n \in J \cap K \cap L$ but $n-1 \not \in J$. Then $Z_J$ is a subset of $\cC \times \{a, a \pm 1\}$ (with the choice of $a+1$ or $a-1$ depending on the color of the chambers in $Z_J$). Indeed, since $n-1 \notin J$, all chambers in $Z_J$ must have the same color. Thus $Z_K$ and $Z_L$ also intersect this same set. Furthermore, from the definition of $r_n$ we see that $Z_J \cap Z_K \cap Z_L = (Z_{\bar{J}} \cap Z_{\bar{K}} \cap Z_{\bar{L}}) \times \{a, a\pm1\}$. Again, $\cC$ is a hypertope, so $Z_{\bar{J}} \cap Z_{\bar{K}} \cap Z_{\bar{L}}$ is a $(\bar{J} \cap \bar{K} \cap \bar{L})$-cell of $\cC$, and since the matching $r_n$ commutes with every other matching $r_i$ with $i < n-1$, it follows that $Z_J \cap Z_K \cap Z_L$ is a $(\{n\} \cup (\bar{J} \cap \bar{K} \cap \bar{L}))$-cell; that is, a $(J \cap K \cap L)$-cell.

    Finally, suppose that $n \in J \cap K \cap L$ and $n-1 \in J$. If $n-1 \notin K$ or $L$, the same argument as above holds symmetrically. Thus, we can assume $\{n-1, n\} \subseteq J \cap K \cap L$. Then $Z_J$ contains chambers of both colors, and it follows that $Z_J = Z_{\bar{J}} \times \mathbb Z_k$, and similarly for $Z_K$ and $Z_L$. Thus $Z_J \cap Z_K \cap Z_L = (Z_{\bar{J}} \cap Z_{\bar{K}} \cap Z_{\bar{L}}) \times \mathbb Z_k$, and the argument is similar to before.
    \end{proof}

    Thus, if $\cH$ is a hypertope, then $\cH|k$, the incidence system associated to $\cC|k$, is also a hypertope. We will occasionally blur the line  
    between these two ideas, mostly by using $\cC|k$ to compute things about $\cH|k$. 
    
    \begin{proposition}
    \label{flat-ext-diagram}
    If $\cH$ is a hypertope with diagram $D$ and $k \geq 4$ is even, 
    then $\cH|k$ has a diagram obtained from $D$ by adding a new node $n$ connected to type $n-1$ with an edge of label $k$.
    \end{proposition}

    \begin{proof}
    It is clear from the construction that the $n$-elements of $\cC|k$ are isomorphic to $\cC$, so we only need to show that the $\{i,n\}$-cells of $\cC|k$ are digons if $i < n-1$ and $k$-gons if $i = n-1$. This is clear from the construction: every $\{i,n\}$-cell has $4$ chambers if $i < n-1$ and $2k$ chambers if $i = n-1$.
    \end{proof}

    \begin{proposition}
    \label{flat-ext-tight}
    If $\cH$ is a hypertope and $k \geq 4$ is even, then 
    $\cH|k$ is $(i,n)$-flat for every $i < n-1$, and it is tight if and only if $\cH$ is.
    \end{proposition}

    \begin{proof}
    The $n$-elements of $\cC|k$ are the $k$ copies of $\cC$, and the $i$-elements correspond to the connected components of $\cC|k$ after we remove all edges of label $i$. We note that, for every chamber $(C, a)$ of $\cC|k$, each $\{n-1,n\}$-cell intersects every copy of $\cC$ (that is, every $n$-element). Hence, for $i < n-1$, since each $i$-element contain $\{n-1,n\}$-cells, every $i$-element intersects every $n$-element, which is to say that $\cC|k$ is $(i,n)$-flat.

    For the second part, the necessity follows from Proposition~\ref{flat-ext-diagram} and Proposition~\ref{tight-residues}(b). Conversely, if $\cC$ is tight, then since $\cC$ has $k$ times as many chambers and has a spanning forest for its diagram with one additional edge labeled $k$, it follows that $\cC|k$ is tight.
    \end{proof}

    Next, let us show that the extensions we build this way are as symmetric as we could hope for.

    \begin{proposition}
    \label{flat-ext-auts}
    $\aut(\cH|k) \cong \aut(\cH) \ltimes D_{k/2}$, where $D_{k/2}$ is the dihedral group of order $k$.
    \end{proposition}

    \begin{proof}
    We will work with $\cC$ and $\cC|k$ directly.
    Let us first show how to extend an automorphism $\varphi$ of $\cC$ to an automorphism of $\cC|k$. Recall that $\cC$ is $(n-1)$-element bipartite; let us fix a bicoloring of $\cC$ (say, using red and blue) and color each copy of $\cC$ in $\cC|k$ by the same bicoloring. 
    For $\varphi \in \aut(\cC)$, we define $\tilde{\varphi}$ by 
    \[ (C, j) \tilde{\varphi} = (C\varphi, j) \]
    if $\varphi$ is color-preserving, and
    \[ (C, j) \tilde{\varphi} = (C \varphi, -j) \]
    if $\varphi$ is color-reversing. It is clear that $\tilde{\varphi}$ commutes with $r_i$ for $i < n$. Next, if $\varphi$ is color-preserving, we have
    \[ (C, j)^n \tilde{\varphi} = (C, j \pm1) \tilde{\varphi} = (C \varphi, j \pm 1) = (C\varphi, j)^n = ((C, j) \tilde{\varphi})^n, \]
    where the third equality follows since the color of $C\varphi$ is the same color as $C$, so $r_n$ acts on the second coordinate of $(C,j)$ and $(C\varphi, j)$ in the same way. On the other hand, if $\varphi$ is color-reversing, we have
    \[ (C, j)^n \tilde{\varphi} = (C, j \pm 1) \tilde{\varphi} = (C \varphi, -j \mp 1) = (C\varphi, -j)^n = ((C, j )\tilde{\varphi})^n,\]
    where the third equality follows since $C\varphi$ has the opposite color of $C$. So every automorphism of $\cC$ induces an automorphism of $\cC|k$ that fixes layer $0$.

    Now, from the structure of $\cC|k$ it is clear that there is an automorphism $x$ that sends each $(C,j)$ to $(C, 1-j)$ and an automorphism $y$ that sends each $(C, j)$ to $(C, -1-j)$. These are involutions whose product composition sends $(C, j)$ to $(C, j+2)$, and so they generate a dihedral group of order $k$. 

    The rest of the argument is essentially the same as that of \cite[Proposition 3.9]{flat-extensions}; every element of $\langle x, y \rangle$ fixes the first coordinate of every pair $(C,j)$, and so it is a normal subgroup, and every automorphism of $\cC|k$ can be written as a product of some $\tilde{\varphi}$ with an element of $\langle x, y \rangle$.    
    \end{proof}

    \begin{corollary}
    \label{flat-ext-reg}
    $\cH|k$ is regular if and only if $\cH$ is.
    \end{corollary}

    \begin{proof}
    The necessity is obvious. Conversely, if $\cH$ is regular, then Proposition~\ref{flat-ext-auts} implies that $\aut(\cH|k)$ acts transitively on chambers, and so $\cH|k$ is regular.
    \end{proof}

    Finally, let us show that the construction can be iterated.

    \begin{proposition}
    \label{flat-ext-still-bipartite}
    Let $k$ be a positive even integer. Then $\cC|k$ is $n$-element bipartite, and for every $j \in \{0, \ldots, n-1\}$, if $\cC$ is $j$-element bipartite then so is $\cC|k$.
    \end{proposition}

    \begin{proof}
    Since $k$ is even and $r_n$ changes the layer of a chamber by $\pm 1$, it is clear that every closed path in $\cC|k$ uses $n$-edges an even number of times. Furthermore, if $\cC$ is $j$-element bipartite, then the induced bipartition on $\cC|k$ (where every layer is bipartitioned in the same way) is a witness that $\cC|k$ is $j$-element bipartite.
    \end{proof}

    Putting everything together, we have proven the following. 

    \begin{theorem}
    \label{flat-ext-thm}
    Suppose $\cH$ is a hypertope with diagram $D$ that is $(n-1)$-element bipartite. Then:
    \begin{enumerate}
        \item $\cH|k$ is a hypertope.
        \item The diagram of $\cH|k$ is obtained from $D$ by adding a new node $n$ connected to type $n-1$ with an edge of label $k$.
        \item $\cH|k$ is tight if and only if $\cH$ is.
        \item $\cH|k$ is regular if and only if $\cH$ is.
        \item $\cH|k$ is $n$-element bipartite and is $j$-element bipartite for each $j$ such that $\cH$ is $j$-element bipartite.
    \end{enumerate}
    \end{theorem}

    \begin{corollary}
    \label{one-odd-tree}
    Suppose $D$ is a tree. Suppose that $D$ has only even labels, or that it has the form
\begin{center}
\begin{tikzpicture}[
    every node/.style={circle, draw, inner sep=0pt, minimum size=5mm},
    thick
]
    \node (A) at (0,0) {};
    \node (B) at (1.5,0) {};
    \node (C) at (3,0) {};

    \draw (A) -- node[draw=none, above, circle] {$p$} (B);
    \draw (B) -- node[draw=none, above, circle] {$q$} (C);

    \draw[dashed] (4,0) ellipse (1.5cm and 0.8cm);
\end{tikzpicture},
\end{center}

    where $p$ is odd and $q \equiv 2$ (mod $4$). 
    Then there is a tight hypertope with diagram $D$. 
    \end{corollary}

    \begin{proof}
    Our goal is to find a facet-bipartite, tight polyhedron to start with. If there is a single odd label $p$ with $q \equiv 2$ (mod $4$), then the tight polyhedron of type $\{q,p\}$ described in \cite[Theorem 5.1]{tight-polytopes} has a bipartite graph, so its dual is a facet-bipartite tight polyhedron of type $\{p,q\}$. Similarly, if every label of $D$ is even and we pick $p$ and $q$ to be arbitrary adjacent labels, then the tight polyhedron of type $\{p,q\}$ described in \cite[Theorem 5.1]{tight-polytopes} is also facet-bipartite. Then by Theorem~\ref{flat-ext-thm}, we can build a tight hypertope with diagram $D$ by repeatedly adding a new type and a new edge with an even label connected to this new type and an existing type.
    \end{proof}

    \begin{remark}
    In \cite[Example 4.1]{cayley-ext}, flat extensions of facet-bipartite polytopes were also described in terms of a voltage graph as follows. Start with a facet-bipartite $n$-polytope (or more precisely, with its chamber system), and color the chambers according to a bipartition of the facets. Add semi-edges of label $n$ to every chamber. Then take a dihedral group $\langle x, y \rangle$ and give voltage $x$ to semi-edges attached to chambers of one color and voltage $y$ to the remaining semi-edges. Then the derived graph is precisely the flat extension described in \cite[Section 3.1]{flat-extensions}.

    We expect that this same construction should work more generally here with hypertopes, but proving this is beyond the scope of this paper.
    \end{remark}
        
\section{Constructing tight hypertopes using twisting}
\label{sec:twisting}

    The construction of Section~\ref{sec:flat-ext} can be applied to hypertopes regardless of any symmetry, and by \cref{flat-ext-thm}, it can build tight regular $(n+1)$-hypertopes from tight regular $n$-hypertopes. In this section, we explore a purely algebraic way to build these regular tight hypertopes.



    In order to do so, we will use a construction called \emph{twisting}. The notion of twisting was originally defined in some very specific cases for polytopes (see \cite{arp}). It was extended to general coset geometries in \cite{AmalgamsHypertopes}, and generalized even further in \cite{WreathCosetGeom}. As we only need a specific instance of twisting in this paper, we will introduce only a simplified version of the construction here. For more details, we advise to consult \cite{WreathCosetGeom}.

    Let $\cH$ be a regular $n$-hypertope. 
    We will need some basic facts about the automorphism groups of regular hypertopes; the full details can be found in \cite{hypertopes}. Choose a base chamber $c$ of $\cH$. As $\cH$ is thin, for each $i = 0,1,\ldots,n-1$, there is a unique chamber $c^i$ that is $i$-adjacent to $c$. Because $\cH$ is regular, there must then exist an automorphism $\alpha_i$ sending $c$ to $c^i$. Also, this automorphism is unique and must be an involution. Moreover, it turns out that $\aut(\cH) = \langle \alpha_0,\alpha_1,\ldots, \alpha_{n-1}\rangle$. Actually, the pair $(\aut(\cH), \{\alpha_0,\alpha_1,\ldots, \alpha_{n-1} \})$ forms a $C$-group,
    but that is not crucial for understanding the construction (for more details on $C$-groups and their relation with regular hypertopes, we refer to~\cite{hypertopes}).

    Let $D_m = \langle x,y\rangle$ be a dihedral group of order $2m$, which is the automorphism group of a regular $m$-gon, and let $\aut(\cH) = \langle \alpha_0,\alpha_1,\ldots, \alpha_{n-1}\rangle$ be the automorphism group of a regular hypertope. Choose any $j \in \{0,1,\ldots, m-1\}$ and let $\varphi_j \colon \aut(\cH) \to \aut(D_m)$ be the map defined as follows:
    \begin{itemize}
        \item If $i \neq j$, we set $\varphi_j(\alpha_i)$ to be the identity.
        \item If $i =j$, we set $\varphi_j(\alpha_j)$ to be the order two automorphism switching $x$ and $y$.
    \end{itemize}
    This map $\varphi_j$ extends to a group homomorphism, and thus to an action of $\aut(\cH)$ on $D_m$, if and only if all relators in $\aut(\cH)$ with respect to the generating set $\{\alpha_0,\alpha_1,\cdots,\alpha_{m-1}\}$ contain an even number of $\alpha_j$'s. Suppose that this is the case and let $G = D_m \rtimes_{\varphi_j} \aut(\cH)$.

    \begin{theorem}\label{thm:twistExtension}
        The group $G = \langle \alpha_0,\alpha_1, \ldots, \alpha_{n-1}, x \rangle$ is the automorphism group of a regular $(n+1)$-hypertope $\cH'$.
        Moreover, the diagram of $\cH'$ is obtained from the diagram of $\cH$ by adding a new node connected by an edge labeled $2m$ to the $j$-node of the diagram of $\cH$. Finally, the hypertope $\cH'$ is tight if $\cH$ is tight.
    \end{theorem}

    \begin{proof}
        Clearly, $G = D_m \rtimes_{\varphi_j} \aut(\cH)$ is generated by $\{\alpha_0,\alpha_1, \ldots, \alpha_{n-1}, x\}$ since $y$ can be recovered by letting $\alpha_j$ act on $x$. For ease of notations, set $\alpha_n = x$. The first part of the claim is proved by showing that the coset geometry $(G,(G_i)_{i= 0,1,\ldots,n})$ is a regular hypertope. This is precisely the content of \cite[Theorem 3.16, Corollary 3.17]{WreathCosetGeom}, in the specific case where $B = \aut(\cH)$ and $A = D_m = \langle x, y\rangle$.

        For the diagrams, notice that the diagram of $\cH'$ is characterized by the orders in $G = D_m \rtimes_{\varphi_j}\aut(\cH)$ of the elements $\alpha_i\alpha_j$ for $i \neq j \in \{0,1, \cdots,n\}$. It is clear that if both $i$ and $j$ are different from $n$, so that $\alpha_i$ and $\alpha_j$ are generators of $\aut(\cH)$, the order of $\alpha_i\alpha_j$ in $G$ is the same as it was in $\aut(\cH)$. Thus, the sub-diagram of $\cH'$ spanned by all nodes except the node $n$ is the diagram of $\cH$. Finally, we need to compute the order of $\alpha_i \alpha_n$ for $i =0,1,\ldots, n-1$. First of all, if $i \neq j$, the generator $\alpha_i$ acts trivially on $D_m$ and it is thus clear that $\alpha_i$ and $\alpha_n$ commute in $G$. If $i =j$, let us properly identify $\alpha_n$ with $(\alpha_n, e)$ and $\alpha_j$ with $(e,\alpha_j)$ in $G = D_m \rtimes_{\varphi_j} \aut(\cH)$. Then, we have that $(e,\alpha_j)(\alpha_n,e)(e,\alpha_j)(\alpha_n,e) = (y,\alpha_j)(e,\alpha_j)(x,e) = (y,e)(x,e) = (yx,e)$. Hence, the order of $\alpha_j\alpha_n$ must be twice the order of $yx$, which is $m$.

        Suppose now that $\cH$ is tight. This means that its diagram is a tree $T$ and that the number of chambers of $\cH$ is equal to twice the product of all labels appearing on the edges of the tree. The diagram of $\cH'$ is also a tree $T'$, with an additional edge. The weight $w(T')$ is thus equal to $w(T) \cdot 2m$ where $2m$ is the label of the new edge. The number of chambers of $\cH'$ is equal to the cardinality of $G = D_m \rtimes_{\varphi_i} \aut(\cH)$ which is equal to $ 2m \cdot 2w(T) = 4 m w(T) =2 w(T')$. This shows that $\cH'$ is tight whenever $\cH$ is tight. 
    \end{proof}
    
We now note that, when starting from a regular hypertope $\cH$, the flat-extension described in Section~\ref{sec:flat-ext} is exactly the construction obtained from Theorem~\ref{thm:twistExtension}. Let us denote by $\cT(D_{m},\cH)$ the hypertope $\cH'$ obtained in Theorem \ref{thm:twistExtension} from $\aut(\cH)$ and $D_m$.

\begin{theorem}\label{thm:twisting_is_flat_ext}
    Let $\cH$ be a regular $n$-hypertope which is $(n-1)$-element bipartite and with diagram $D$.
    Then $\cT(D_{m},\cH) = \cH|2m$.
    Moreover, if the presentation of $\aut(\cH)$ is $\langle \alpha_0,\ldots,\alpha_{n-1}\mid R_\alpha\rangle$, then $\aut(\cH|2m)$ has presentation $$\langle \alpha_0,\ldots,\alpha_{n-1},x\mid R_\alpha, x^2=(xa_{n-1})^{2m}=[x,\alpha_i]= [[\alpha_{n-1},x],\alpha_i]=1,\textnormal{ for }i<n-1\rangle$$ 
\end{theorem}
\begin{proof}
    Let $\cH$ be a regular $n$-hypertope $\cH$ which is $(n-1)$-element bipartite and let $\cC$ be its chamber system. Let $\aut(\cH)=\langle \alpha_0,\ldots,\alpha_{n-1}\rangle$, as before.
    Since $\cC$ is $(n-1)$-element bipartite, all closed loops starting in a chamber $c\in \cC$ must cross an even number of $(n-1)$-labeled edges. As closed loops in $\cC$ correspond to words on the alphabet $\{\alpha_0,\ldots,\alpha_{n-1}\}$ which are sent to the identity in $\aut(\cH)$, this is equivalent to all relators of $\aut(\cH)$ containing an even number of $\alpha_{n-1}$. Hence, we see that the condition for $\varphi_{n-1}$ to extend to a group homomorphism is equivalent to $\cH$ being $(n-1)$-element bipartite.

    Let $D_m = \langle x,y\mid x^2=y^2=(xy)^m=1\rangle$ as before, and let $G = \aut(\cT(D_m,\cH))$. Since $G = D_m \rtimes \aut(\cH)$, and $y = \alpha_{n-1} x \alpha_{n-1}$, we can eliminate $y$ from a presentation for $G$ to get the presentation
    $$\langle S_\alpha,x\mid R_\alpha, x^2=(x\alpha_{n-1})^{2m} =1, (x\alpha_{i})^2 = (\alpha_{n-1}x\alpha_{n-1}\alpha_{i})^2=1,\textnormal{ for } i<n-1\rangle,$$
    where $S_\alpha = \{\alpha_0,\ldots,\alpha_{n-1}\}$ and $R_\alpha$ is the set of relators of $\aut(\cH)$ with respect to $S_\alpha$.
    Note that $(\alpha_{n-1}x\alpha_{n-1}\alpha_{i})^2=1$ can be alternatively written as $[[\alpha_{n-1},x],\alpha_i]=1$ and $(x\alpha_i)^2=1$ as $[x,\alpha_i]=1$.

    Now consider the regular hypertope obtained from the construction $\cH|2m$, with automorphism group $\aut(\cH|2m) = \aut(\cH)\ltimes D_m$. As described in Proposition~\ref{flat-ext-auts}, $\aut(\cH|2m) = \langle S_\alpha, x, y\rangle$, where $x$ is defined as the automorphism $(C,j)\mapsto (C,1-j)$ and $y$ as the automorphism $(C,j)\mapsto(C,-1-j)$. 
    We will now prove that $\aut(\cH|2m)$ satisfies all the relations of $G$. First, note that all generators and relations of $\aut(\cH)$ are still true in $\aut(\cH|2m)$. Moreover, one can easily check that $y = \alpha_{n-1}x\alpha_{n-1} = \varphi_{n-1}(\alpha_{n-1})(x)$, meaning that $(x\alpha_{n-1}x\alpha_{n-1})^{m} = (x\alpha_{n-1})^{2m}$ must fix every chamber. Additionally, for all $\alpha_i$ with $i<n-1$, we have that $(x\alpha_i)^2$ and $(y\alpha_i)^2$ fixes every chamber.
    Substituting $y=\alpha_{n-1}x\alpha_{n-1}$, we get that $(\alpha_{n-1}x\alpha_{n-1}\alpha_i)^2$ fixes every chamber.
    Hence $\aut(\cH|2m)$ is a quotient of the group $G$. However, since both are tight, they must have the same order, and hence, are the same group. Therefore, $\cH|2m$ is the same tight regular hypertope as $\cT(D_m,\cH)$.
\end{proof}

Starting from a coset geometry of $k$-gon, with $k$ even, we can iterate the construction above to create a regular tight hypertope whose diagram is any tree with all even labels.

\begin{corollary}\label{cor:presentation}
    Let $T$ be any labeled tree, with all labels equal to even integers. Then there exists a regular tight hypertope whose diagram is $T$ and whose automorphism group has presentation
\[
\left\langle
 \alpha_i, i\in V(T)
 \ \middle|\
\begin{array}{l}
\alpha_i^{2}= 1,
[ \alpha_i, \alpha_j] = 1\textnormal{ for }\{i,j\}\notin E(T),\\[2mm]
(\alpha_i \alpha_j)^{p_{ij}} = 1\textnormal{ for }\{i,j\}\in E(T)\textnormal{ with even label } p_{ij},\\[2mm]
[[ \alpha_i, \alpha_j], \alpha_k] = 1\textnormal{ for all distinct }i,j,k \in V(T)
\end{array}
\right\rangle.
\]
\end{corollary}

\begin{proof}
    
     Let $k$ be an even number and $D_k = \langle \alpha_0,\alpha_1 \rangle$ be a dihedral group of order $2k$. Then $D_k$ is the automorphism group of a regular $k$-gon, which is a tight regular polytope, and thus hypertope. Notice that, since $k$ is even, every defining relation 
     of $D_k$ contains an even number of $\alpha_0$ and an even number of $\alpha_1$. We can thus let $D_k$ act on any other dihedral group $D_m = \langle x,y \rangle$ via $\varphi_0$ or $\varphi_1$. Doing so, by Theorem \ref{thm:twistExtension} we form a regular tight hypertope $\cH$ (actually a polytope in this case) whose diagram is a path of length two with labels $k$ and $2m$. Since $k$ and $m$ can be chosen arbitrarily, the theorem is proven for $|T| = 3$.
    As mentioned in Theorem~\ref{thm:twistExtension}, the choice of $\varphi_i$ adds a new edge of label $2m$ to the previous diagram.
    
     The automorphism group of the hypertope $\cH$ is $D_m\rtimes_{\varphi_i} D_{k} = \langle \alpha_0, \alpha_1, \alpha_2 = x \rangle$, with presentation
    $$\langle \alpha_0,\alpha_1,\alpha_2\mid \alpha_l^2=(\alpha_0\alpha_1)^k=(\alpha_1\alpha_2)^{2m}=(\alpha_0\alpha_2)^2= [[\alpha_1,\alpha_2],\alpha_0]=1\rangle, $$
    if $\varphi_i=\varphi_1$, otherwise
     $$\langle \alpha_0,\alpha_1,\alpha_2\mid \alpha_l^2=(\alpha_0\alpha_1)^k=(\alpha_0\alpha_2)^{2m}=(\alpha_1\alpha_2)^2= [[\alpha_0,\alpha_2],\alpha_1]=1\rangle. $$
     The relators for this group all contain an even number of any of the three generators since the relations introduced by semi-direct products are all of the form $zwz^{-1} = \varphi_i(z)(w)$.
     Therefore, we can now use $\cH$ as our starting hypertope, let it act on any dihedral group and form a new regular tight hypertope whose diagram has one more edge with even label. 
     Iterating this process, we can construct tight regular hypertopes whose diagram is any tree whose labels are all even integers. Additionally, by Theorems~\ref{thm:twistExtension} and~\ref{thm:twisting_is_flat_ext}, this iteration results in the presentation given in the statement of the corollary. (Note that, if the subdiagram with vertices $\{i,j,k\}$ is not connected, then in fact the relation $[[\alpha_i, \alpha_j], \alpha_k] = 1$ is already implied by the other relations.)
\end{proof}
    
Note that if we are given a regular tight hypertope whose diagram has some odd labels, we can still attempt to use this process to reach bigger trees. We will only need to check that $\varphi_i$ is a homomorphism. Notice that if the vertex $i$ of the diagram is the end point of any odd labeled edge, $\varphi_i$ will never be a homomorphism. So we can only hope to glue more edges to vertices that only have even labeled edges to start with.

Finally, we show that whenever all labels are equal to $4$, there is a unique tight regular hypertope, just as it is the case for polytopes.

    
\begin{theorem}
    There is a unique regular tight hypertope whose diagram is a tree with all labels equal to 4. 
\end{theorem}
\begin{proof}
    
    Let $T$ be a tree whose edges are all labeled by $4$. By Corollary \ref{cor:presentation}, we know that there exists a tight regular hypertope with diagram $T$ and that its automorphism group $G$ has presentation
    \[
     G :=
    \left\langle
    \rho_i , \ i\in V(T)
    \ \middle|\
    \begin{array}{ll}
    \rho_i^2 = 1
    & \text{for all } i\in I,\\[2mm]
    [[\rho_i,\rho_j],\rho_k]=1
    & \text{for all } i,j,k\in V(T),\\[2mm]
    [\rho_i,\rho_j]=1
    & \text{whenever } \{i,j\}\notin E(T)
    \end{array}
    \right\rangle .
    \]
    Notice that this group is precisely the one constructed in \cite{HypertopesDiagramSymmetries} when the underlying graph is a tree.
    
    We now show that the automorphism group of any tight regular hypertope $\cH$ whose diagram is a tree $T$ all of whose labels are $4$ must be isomorphic to $G$. It is clear that $\aut(\cH)$ is generated by $|V|$ involutions $\alpha_1,\alpha_2,\ldots, \alpha_{|V|}$. Moreover, whenever $\alpha_i$ and $\alpha_j$ correspond to non-adjacent vertices of $T$, it must be that $\alpha_i$ and $\alpha_j$ commute. We now show that $[[\alpha_i,\alpha_j],\alpha_k] = 1$ for any $i,j,k \in I$. The only non trivial case to check is when $i-j-k$ forms a path in $T$. Consider the $\{i,j,k\}$-residue of $\cH$. This residue is a thin, residually connected and regular geometry whose diagram is a path of length three with labels equal to four. It must hence be the unique tight polytope $\cP$ of type $\{4,4\}$, which is commonly denoted by $\{4,4\}_{(2,0)}$. (See, for example, the data at \cite{atlas}.) The automorphism group of this polytope is precisely the group $G$ when $T$ is the only tree with three vertices (this can be verified in a computer algebra system). In particular, the desired commutator relation holds in $\cP$. Since $\cH$ is a regular hypertope, this residue can be identified with the coset geometry $(\langle\alpha_i,\alpha_j,\alpha_k\rangle, (\langle \alpha_j,\alpha_k\rangle,\langle \alpha_i,\alpha_k\rangle ,\langle \alpha_i,\alpha_j\rangle))$. Therefore, we conclude that $[[\alpha_i,\alpha_j],\alpha_k] = 1$.

    We have shown that the automorphism group of any tight regular hypertope $\cH$ with diagram $T$ is generated by $|V|$ involutions and that these involutions satisfy all the relations in the presentation of $G$. Hence, the group of $\cH$ must be a quotient of $G$. But since $\cH$ is tight and regular, $|\aut(\cH)| = 2 \cdot 4^{|E|} = |G|$. Therefore, we can indeed deduce that $\aut(\cH) \cong G$.
\end{proof}

\section*{Acknowledgments}

The authors thank Maria Elisa Fernandes for helpful discussions about small hypertopes and $C$-groups. 

The second author is funded by an Action de Recherche Concert\'{e}e -- ARC -- from the Communaut\'{e}e Fran\c{c}aise Wallonie-Bruxelles. 

The third author is a Postdoctoral Researcher of the Fonds de la Recherche Scientifique -- FNRS.

\bibliographystyle{plain}
\bibliography{flat}

@article {tight3,
		AUTHOR = {Cunningham, Gabe and Pellicer, Daniel},
		 TITLE = {Classification of tight regular polyhedra},
	   JOURNAL = {J. Algebraic Combin.},
	  FJOURNAL = {Journal of Algebraic Combinatorics. An International Journal},
		VOLUME = {43},
		  YEAR = {2016},
		NUMBER = {3},
		 PAGES = {665--691},
		  ISSN = {0925-9899},
	   MRCLASS = {52B15 (05E18)},
	  MRNUMBER = {3482443},
	MRREVIEWER = {Peter A. Brooksbank},
		   DOI = {10.1007/s10801-015-0649-3},
		   URL = {http://dx.doi.org/10.1007/s10801-015-0649-3},
	}

@article {tight-polytopes,
    AUTHOR = {Cunningham, Gabe},
     TITLE = {Minimal equivelar polytopes},
   JOURNAL = {Ars Math. Contemp.},
  FJOURNAL = {Ars Mathematica Contemporanea},
    VOLUME = {7},
      YEAR = {2014},
    NUMBER = {2},
     PAGES = {299--315},
      ISSN = {1855-3966},
   MRCLASS = {52B05 (05C25 51M20 52B15)},
  MRNUMBER = {3240439},
MRREVIEWER = {Thomas Marco Connor},
       DOI = {10.26493/1855-3974.357.422},
       URL = {https://doi.org/10.26493/1855-3974.357.422},
}

@article {flat-extensions,
    AUTHOR = {Cunningham, Gabe},
     TITLE = {Flat extensions of abstract polytopes},
   JOURNAL = {Art Discrete Appl. Math.},
  FJOURNAL = {The Art of Discrete and Applied Mathematics},
    VOLUME = {4},
      YEAR = {2021},
    NUMBER = {3},
     PAGES = {Paper No. 3.06, 14},
   MRCLASS = {52B05 (52B15)},
  MRNUMBER = {4312690},
       DOI = {10.26493/2590-9770.1351.f7e},
       URL = {https://doi-org.ezproxy.lib.umb.edu/10.26493/2590-9770.1351.f7e},
}

@article{cayley-ext,
title = {Cayley extensions of maniplexes and polytopes},
journal = {Journal of Combinatorial Theory, Series A},
volume = {212},
pages = {106000},
year = {2025},
issn = {0097-3165},
doi = {https://doi.org/10.1016/j.jcta.2024.106000},
url = {https://www.sciencedirect.com/science/article/pii/S0097316524001390},
author = {Gabe Cunningham and Elías Mochán and Antonio Montero}
}

@article {tight2,
		AUTHOR = {Conder, Marston and Cunningham, Gabe},
		 TITLE = {Tight orientably-regular polytopes},
	   JOURNAL = {Ars Math. Contemp.},
	  FJOURNAL = {Ars Mathematica Contemporanea},
		VOLUME = {8},
		  YEAR = {2015},
		NUMBER = {1},
		 PAGES = {68--81},
		  ISSN = {1855-3966},
	   MRCLASS = {52B15 (05E18 52B05)},
	  MRNUMBER = {3281121},
	MRREVIEWER = {Daniel Pellicer},
	}

@article{tight-chiral-polyhedra,
	  title={Tight chiral polyhedra},
	  author={Cunningham, Gabe},
	  journal={Combinatorica},
	  pages={115--142},
	  year={2018},
	  volume = {38},
	  number = {1},
	  publisher={Springer Berlin Heidelberg},
	  doi={10.1007/s00493-016-3505-8},
}

@book {arp,
		AUTHOR = {McMullen, Peter and Schulte, Egon},
		 TITLE = {Abstract regular polytopes},
		SERIES = {Encyclopedia of Mathematics and its Applications},
		VOLUME = {92},
	 PUBLISHER = {Cambridge University Press},
	   ADDRESS = {Cambridge},
		  YEAR = {2002},
		 PAGES = {xiv+551},
		  ISBN = {0-521-81496-0},
	   MRCLASS = {52B15 (20F55 51F15 51M20)},
	  MRNUMBER = {1965665 (2004a:52020)},
	MRREVIEWER = {Rade {\v{Z}}ivaljevi{\'c}},
		   DOI = {10.1017/CBO9780511546686},
		   URL = {http://doi.org/10.1017/CBO9780511546686},
	}

@article {tight-chiral-polytopes,
    AUTHOR = {Cunningham, Gabe and Pellicer, Daniel},
     TITLE = {Tight chiral polytopes},
   JOURNAL = {J. Algebraic Combin.},
  FJOURNAL = {Journal of Algebraic Combinatorics. An International Journal},
    VOLUME = {54},
      YEAR = {2021},
    NUMBER = {3},
     PAGES = {837--878},
      ISSN = {0925-9899},
   MRCLASS = {Prelim},
  MRNUMBER = {4330415},
       DOI = {10.1007/s10801-021-01023-z},
       URL = {https://doi-org.ezproxy.lib.umb.edu/10.1007/s10801-021-01023-z},
}

@book{buekenhout,
  title={Diagram geometry: related to classical groups and buildings},
  author={Buekenhout, Francis and Cohen, Arjeh M},
  volume={57},
  year={2013},
  publisher={Springer Science \& Business Media}
}

@article{AmalgamsHypertopes,
title = {From group operations to geometric structures: Amalgamations, HNN-extensions, and twisting in coset geometries},
journal = {Journal of Algebra},
volume = {707},
pages = {83-144},
year = {2026},
issn = {0021-8693},
doi = {https://doi.org/10.1016/j.jalgebra.2026.05.032},
url = {https://www.sciencedirect.com/science/article/pii/S0021869326003339},
author = {Claudio Alexandre Piedade and Philippe Tranchida},
}

@article{WreathCosetGeom,
  title={The geometry of wreath and semi-direct products},
  author={Piedade, Claudio Alexandre and Tranchida, Philippe},
  journal={arXiv preprint arXiv:2605.04646},
  year={2026}
}

@incollection {Tits1982,
    AUTHOR = {Tits, J.},
     TITLE = {A local approach to buildings},
 BOOKTITLE = {The geometric vein},
     PAGES = {519--547},
 PUBLISHER = {Springer, New York-Berlin},
      YEAR = {1981},
   MRCLASS = {51D15 (20G99 51J05)},
  MRNUMBER = {661801},
MRREVIEWER = {Mark A. Ronan},
}

@incollection {Tits1957,
    AUTHOR = {Tits, J.},
     TITLE = {Sur les analogues alg\'{e}briques des groupes semi-simples
              complexes},
 BOOKTITLE = {Colloque d'alg\`ebre sup\'{e}rieure, tenu \`a {B}ruxelles du 19 au 22
              d\'{e}cembre 1956},
    SERIES = {Centre Belge de Recherches Math\'{e}matiques},
     PAGES = {261--289},
 PUBLISHER = {\'{E}tablissements Ceuterick, Louvain},
      YEAR = {1957},
   MRCLASS = {50.00 (22.00)},
  MRNUMBER = {108765},
MRREVIEWER = {H. Freudenthal},
}

@article {hypertopes,
    AUTHOR = {Fernandes, Maria Elisa and Leemans, Dimitri and Weiss, Asia
              Ivi\'c},
     TITLE = {Highly symmetric hypertopes},
   JOURNAL = {Aequationes Math.},
  FJOURNAL = {Aequationes Mathematicae},
    VOLUME = {90},
      YEAR = {2016},
    NUMBER = {5},
     PAGES = {1045--1067},
      ISSN = {0001-9054,1420-8903},
   MRCLASS = {51E24 (20F05 52B11)},
  MRNUMBER = {3547707},
MRREVIEWER = {Hendrik\ Van Maldeghem},
       DOI = {10.1007/s00010-016-0431-1},
       URL = {https://doi.org/10.1007/s00010-016-0431-1},
}

@misc{HypertopesDiagramSymmetries,
      title={Hypertopes with prescribed diagram symmetries}, 
      author={Dimitri Leemans and Pablo Spiga and Philippe Tranchida},
      year={2026},
      eprint={2610.09711},
      archivePrefix={arXiv},
      primaryClass={math.GR},
      url={https://arxiv.org/abs/2610.09711}, 
}

@article{leemans2025groups,
  title={Groups represented by incidence geometries},
  author={Leemans, Dimitri and Stokes, Klara and Tranchida, Philippe},
  journal={arXiv preprint arXiv:2506.10464},
  year={2025}
}

@article {atlas,
	   author = {Hartley, Michael I.},
	   affiliation = {University of Nottingham (Malaysia Campus) Jalan Broga, Semenyih, 43500, Malaysia Jalan Broga, Semenyih, 43500, Malaysia},
	   title = {An Atlas of Small Regular Abstract Polytopes},
	   journal = {Periodica Mathematica Hungarica},
	   publisher = {Akadémiai Kiadó, co-published with Springer Science+Business Media B.V., Formerly Kluwer Academic Publishers B.V.},
	   issn = {0031-5303},
	   pages = {149--156},
	   volume = {53},
	   issue = {1},
	   url = {http://dx.doi.org/10.1007/s10998-006-0028-x},
	   note = {available online at \url{http://www.abstract-polytopes.com/atlas/}},
	   year = {2006}
	}

\end{document}